\documentclass[a4paper]{amsart}  

\usepackage{amssymb} \usepackage{amscd} \usepackage{times}

\def\zbb{\mathbb{Z}}  
  
  \def\phi{\varphi}
 \def\p1{{\mathbb{P}^1_\zbb}}

\linethickness{0.2mm}

\usepackage{tikz} 

\begin{document}

\title{ On the mass of the exterior blow-up points.}

\author{Samy Skander Bahoura}

\address{Departement de Mathematiques, Universite Pierre et Marie Curie, 2 place Jussieu, 75005, Paris, France.}
              
\email{samybahoura@yahoo.fr} 

\date{}

\maketitle

\begin{abstract}

We consider the following problem on open set  $ \Omega $ of $ {\mathbb R}^2 $:

$$  \left \{ \begin {split} 
      -\Delta u_i & = V_i e^{u_i} \,\, &\text{in} \,\, &\Omega \subset {\mathbb R}^2, \\
                  u_i  & = 0  \,\,             & \text{in} \,\,    &\partial \Omega.              
\end {split}\right.
$$
We assume that :

$$ \int_{\Omega} e^{u_i} dy  \leq C, $$

and,

$$ 0 \leq V_i \leq b  < + \infty $$

On the other hand, if  we assume that $ V_i  $ $ s- $holderian with $ 1/2< s \leq 1$, then, each exterior blow-up point is simple. As application, we have a compactness result for the case when:

$$ \int_{\Omega}V_i e^{u_i} dy \leq 40\pi-\epsilon , \,\, \epsilon >0 $$

\end{abstract}

\section{Introduction and Main Results} 

We set $ \Delta = \partial_{11} + \partial_{22} $  on open set $ \Omega $ of $ {\mathbb R}^2 $ with a smooth boundary.

\bigskip

We consider the following problem on $ \Omega \subset {\mathbb R}^2 $:

$$ (P) \left \{ \begin {split} 
      -\Delta u_i & = V_i e^{u_i} \,\, &\text{in} \,\, &\Omega \subset {\mathbb R}^2, \\
                  u_i  & = 0  \,\,             & \text{in} \,\,    &\partial \Omega.              
\end {split}\right.
$$

We assume that,

$$ \int_{\Omega} e^{u_i} dy  \leq C, $$

and,

$$ 0 \leq V_i \leq b  < + \infty $$

The previous equation is called, the Prescribed Scalar Curvature equation, in
relation with conformal change of metrics. The function $ V_i $ is the
prescribed curvature.

\bigskip

Here, we try to find some a priori estimates for sequences of the
previous problem.

\smallskip

Equations of this type were studied by many authors, see [5-8, 10-15]. We can
see in [5], different results for the solutions of those type of
equations with or without boundaries conditions and, with minimal
conditions on $ V $, for example we suppose $ V_i \geq 0 $ and  $ V_i
\in L^p(\Omega) $ or $ V_ie^{u_i} \in L^p(\Omega) $ with $ p \in [1,
+\infty] $. 

Among other results, we  can see in [5], the following important Theorem,

\smallskip

{\bf Theorem A} {\it (Brezis-Merle [5])}.{\it If $ (u_i)_i $ and $ (V_i)_i $ are two sequences of functions relatively to the previous problem $ (P) $ with, $ 0 < a \leq V_i \leq b < + \infty $, then, for all compact set $ K $ of $ \Omega $,

$$ \sup_K u_i \leq c = c(a, b, m, K, \Omega) \,\,\, {\rm if } \,\,\, \inf_{\Omega} u_i \geq m. $$}

A simple consequence of this theorem is that, if we assume $ u_i = 0 $ on $ \partial \Omega $ then, the sequence $ (u_i)_i $ is locally uniformly bounded. We can find in [5] an interior estimate if we assume $ a=0 $, but we need an assumption on the integral of $ e^{u_i} $, precisely, we have in [5]:

\smallskip

{\bf Theorem B} {\it (Brezis-Merle [5])}.{\it If $ (u_i)_i $ and $ (V_i)_i $ are two sequences of functions relatively to the previous problem $ (P) $ with, $ 0 \leq V_i \leq b < + \infty $, and,

$$ \int_{\Omega} e^{u_i} dy  \leq C, $$

then, for all compact set $ K $ of $ \Omega $,

$$ \sup_K u_i \leq c = c(b, C, K, \Omega). $$}

\bigskip

If, we assume $ V $ with more regularity, we can have another type of estimates, $ \sup + \inf $. It was proved, by Shafrir, see [13], that, if $ (u_i)_i, (V_i)_i $ are two sequences of functions solutions of the previous equation without assumption on the boundary and, $ 0 < a \leq V_i \leq b < + \infty $, then we have the following interior estimate:

$$ C\left (\dfrac{a}{b} \right ) \sup_K u_i + \inf_{\Omega} u_i \leq c=c(a, b, K, \Omega). $$

We can see in [7], an explicit value of $ C\left (\dfrac{a}{b}\right ) =\sqrt {\dfrac{a}{b}} $. In his proof, Shafrir has used the Stokes formula and an isoperimetric inequality, see [3]. For Chen-Lin, they have used the blow-up analysis combined with some geometric type inequality for the integral curvature.

\bigskip

Now, if we suppose $ (V_i)_i $ uniformly Lipschitzian with $ A $ the
Lipschitz constant, then, $ C(a/b)=1 $ and $ c=c(a, b, A, K, \Omega)
$, see Br\'ezis-Li-Shafrir [4]. This result was extended for
H\"olderian sequences $ (V_i)_i $ by Chen-Lin, see  [7]. Also, we
can see in [10], an extension of the Brezis-Li-Shafrir to compact
Riemann surface without boundary. We can see in [11] explicit form,
($ 8 \pi m, m\in {\mathbb N}^* $ exactly), for the numbers in front of
the Dirac masses, when the solutions blow-up. Here, the notion of isolated blow-up point is used. Also, we can see in [14] refined estimates near the isolated blow-up points and the  bubbling behavior  of the blow-up sequences.

We have in [15]:

\bigskip

{\bf Theorem C} {\it (Wolansky.G.[15])}. {\it If $ (u_i) $ and $ (V_i) $ are two sequences of  functions solutions of the problem $ (P) $  without the boundary condition, with, 

$$ 0 \leq V_i \leq b < + \infty, $$

$$ ||\nabla V_i||_{L^{\infty}(\Omega)}\leq C_1, $$

$$ \int_{\Omega} e^{u_i} dy  \leq C_2, $$

and,

$$ \sup_{\partial \Omega} u_i - \inf_{\partial \Omega } u_i \leq C_3, $$

 the last condition replace the boundary condition.
 
 \smallskip
 
We assume that $ (iii) $ holds in the theorem 3 of [5], then, in the sense of the distributions:

$$  V_i e^{u_i}  \to \sum_{j=0}^m 8 \pi  \delta_{x_j}. $$

in other words, we have:

 $$ \alpha_j = 8 \pi, \,\,\, j=0 \ldots m,  $$

in $ (iii) $ of the theorem 3 of [5].}

\bigskip

To understand the notations, it is interessant to take a look to a previous prints on arXiv, see [1] and [2].

\bigskip

Our main results are:

\bigskip

{\bf Theorem 1}. {\it Assume that, $ V_i  $ is uniformly $ s- $holderian with $ 1/2 < s \leq 1$, and that :

$$ \max_{\Omega} u_i \to + \infty. $$

Then, each exterior blow-up point is simple. 

\bigskip

There are $ m $ blow-ups points on the boundary (perhaps the same) such that:

 $$ \int_{B(x_i^j, \delta_i^j \epsilon')} V_i(x_i^j+ \delta_i^j y) e^{u_i} \to 8 \pi. $$

and,

$$  \int_{\Omega} V_i e^{u_i} \to  \int_{\Omega} V e^{u} + \sum_{j=1}^m 8 \pi \delta_{x_j} . $$
}

and,

\bigskip

{\bf Theorem 2}. {\it Assume that, $ V_i  $ is uniformly $ s- $holderian with $ 1/2 < s \leq 1$, and,

$$ \int_{B_1(0)} V_i e^{u_i}  dy  \leq 40\pi-\epsilon, \,\, \epsilon>0,$$

then we have:

$$ \sup_{\Omega} u_i \leq c=c(b, C, A, s, \Omega). $$

where $ A $ is the holderian constant of $ V_i $.}

\section{Proof of the result:} 

\underbar {Proof of the theorem 1:}

\bigskip

Let's consider the following function on the ball of center $ 0 $ and radius $ 1/2 $; And let us consider $ \epsilon >0 $

$$ v_i(y)=u_i(x_i+\delta_i y)+2\log \delta_i, \quad y \in B(0,1/2) $$

This function is solution of the following equation:

$$ -\Delta v_i= V_i(x_i+\delta_i y) e^{v_i}, \quad y \in B(0,1/2) $$

The function $ v_i $ satisfy the following inequality (without loss of generality):

$$ \sup_{\partial B(0, 1/4)} v_i-\inf_{\partial B(0, 1/4)} v_i \leq C, $$
 
 Let us consider the following  functions: 
 
$$  \left \{ \begin {split} 
       -\Delta v_0^i &= 0 \,\, && \text{in} \,\,B(0,1/4) \\
               v_0^i  &= u_i(x_i+\delta_i y) \,\, && \text{on} \,\, \partial B(0, 1/4).              
\end {split}\right.
$$

By the elliptic estimates we have:

$$ v_0^i \in C^2 (\bar {B}(0, 1/4)). $$

We can write:

$$ -\Delta (v_i-v_0^i)= V_i(x_i+\delta_i y)e^{v_0^i}e^{v_i-v_0^i}= K_1K_2 e^{v_i-v_0^i}, $$

With this notations, we have:

$$ ||\nabla( v_i- v_0^i)||_{L^q(B(0,\epsilon))} \leq C_q. $$

$$ v_i-v_0^i \to G \,\, \text{in} \,\, W_0^{1, q}, $$

And, because, for  $ \epsilon >0 $ small enough:

$$ ||\nabla G||_{L^q(B(0,\epsilon))} \leq \epsilon' << 1, $$

We have, for $ \epsilon >0 $ small enough:

$$ ||\nabla( v_i- v_0^i)||_{L^q(B(0,\epsilon))} \leq 2 \epsilon' << 1. $$

and,

$$ ||\nabla v_i||_{L^q(B(0,\epsilon))} \leq 3 \epsilon' << 1. $$

Set,

$$ u= v_i-v_0^i, \,\, z_1=0, $$

Then, 

$$ -\Delta u= K_1K_2 e^u, \,\, \text{in} \,\, B(0, 1/4), $$

and,

$$ osc (u) = 0. $$





 We use Woalnsky's theorem, see [15]. In fact $ K_2 $ is a $ C^1 $ function uniformly bounded and $ K_1 $ is s-holderian with $  1/2< s \leq 1 $. Because we take the logarithm in $ K $, the part which contain  $ K_2 $ have similar proof  as in this paper we use the Stokes formula. Only the case of $ K_1 $ s-holderian is difficult. For this and without loss of generality, we can assume the $ K=K_1=V_i(x_i+\delta_i y) $. We set:

$$ \Delta \tilde u=  \Delta v_i= \rho= -K e^{\tilde u}=-K_1e^{v_i} $$

Let us consider the following term of Wolansky computations:

$$ \int_{B^{\epsilon}} div ((z-z_1)\rho) \log K+ \int_{\partial B^{\epsilon}} (<(z-z_1) |\nu >\rho) \log K, $$

First, we write:

$$ \int_{B^{\epsilon}} div ((z-z_1)\rho) \log K=2\int_{B^{\epsilon}} \rho \log K +\int_{B^{\epsilon}}  <(z-z_1)| \nabla \rho) \log K $$ 

which we can write as:

$$ -\int_{B^{\epsilon}} div ((z-z_1)\rho) \log K=2\int_{B^{\epsilon}}  K \log K e^u +\int_{B^{\epsilon}}  <(z-z_1)| \nabla u> K \log K e^u + \int_{B^{\epsilon}}  <(z-z_1) |(\nabla K) \log K> e ^u, $$

We can write:

$$ \nabla (K(\log K)-K)=(\nabla K)(\log K) $$

Thus,  and by integration by part we have:

$$ \int_{B^{\epsilon}}  <(z-z_1)|(\nabla K) \log K> e ^u = \int_{B^{\epsilon}}  <(z-z_1) |(\nabla( K \log K-K))> e ^u = $$

$$ = \int_{\partial B^{\epsilon}}  <(z-z_1)|\nu >(K \log K - K )e ^u - 2 \int_{B^{\epsilon}}  ( K \log K-K)e ^u-\int_{B^{\epsilon}}  <(z-z_1)| \nabla u >( K \log K-K) e ^u $$

Thus,

$$ -(\int_{B^{\epsilon}} div ((z-z_1)\rho) \log K+ \int_{\partial B^{\epsilon}} (<(z-z_1) |\nu >\rho) \log K) = $$

$$ = -\int_{\partial B^{\epsilon}}  <(z-z_1)|\nu > K e ^u + \int_{B^{\epsilon}}  <(z-z_1)| \nabla u >K e ^u+2 \int_{B^{\epsilon}}  K  e^u$$

But, we can write the following,

$$ \int_{B^{\epsilon}}  <(z-z_1) |\nabla u >K e ^u =  \int_{B^{\epsilon}}  <(z-z_1) |\nabla u >(K-K(z_1)) e ^u+K(z_1)  \int_{B^{\epsilon}}  <(z-z_1) |\nabla u > e ^u, $$

and, after integration by parts:

$$ K(z_1)  \int_{B^{\epsilon}}  <(z-z_1)| \nabla u > e ^u = K(z_1)\int_{\partial B^{\epsilon}}  <(z-z_1)|\nu >  e ^u- 2K(z_1)\int_{ B^{\epsilon}} e ^u, $$

Finaly, we have, for the Wolansky term:

$$ \int_{B^{\epsilon}} div ((z-z_1)\rho) \log K+ \int_{\partial B^{\epsilon}} (<(z-z_1) |\nu >\rho) \log K = $$

$$ = \int_{B^{\epsilon}}  <(z-z_1)| \nabla u >(K-K(z_1)) e ^u+\left ( 2 \int_{B^{\epsilon}} ( K -K(z_1)) e^u \right )+ $$

$$ + \left ( \int_{\partial B^{\epsilon}}  <(z-z_1)| \nu >(K(z_1)-K) e ^u \right )$$

But, we have soon that if $ K $ is $ s-$holderian with $ 1 \geq s>1/2 $, around each exteriror blow-up we have, the following estimate:

$$\int_{B^{\epsilon}}  <(z-z_1) |\nabla u >(K-K(z_1)) e ^u= $$

$$ =\int_{B(0, \epsilon)}  <(y-z_1) |\nabla v_i >(V_i(x_i+\delta_iy)-V_i(x_i)) e^{v_i} dy= $$ 
$$=\int_{B(x_i, \delta_i\epsilon)}  <(x-x_i) |\nabla u_i >(V_i(x)-V_i(x_i)) e^{u_i} dy=o(1)M_{\epsilon} $$
$$ =  o(1)\int_{B(x_i, \delta_i \epsilon)} V_i e^{u_i} = o(1)\int_{B^{\epsilon}} K e^u, $$

Thus,

$$ \int_{B^{\epsilon}} div ((z-z_1)\rho) \log K+ \int_{\partial B^{\epsilon}} (<(z-z_1)| \nu >\rho) \log K =o(1)M_{\epsilon}= o(1)\int_{B^{\epsilon}} K e^u $$

We argue by contradiction and we suppose that we have around the exterior blow-up point 2 or 3 blow-up points, for example. We prove, as in a previous paper, that, the last quantity tends to 0. But according to Wolansky paper, see [15]:

$$ \int_{B^{\epsilon}} V_i(x_i+ \delta_i y) e^{v_i} \to 8 \pi. $$

Around  each exterior blow-up points, there is one blow-up point.

Consider the following quantity:

$$ B_i=\int_{B(x_i, \delta_i\epsilon)}  <(x-x_i) |\nabla u_i >(V_i(x)-V_i(x_i)) e^{u_i} dy. $$

Suppose that, we have $ m >0 $ interior blow-up points. Consider the blow-up point $ t_i^k $ and the associed set $ \Omega_k $ defined as  the set of the points nearest $ t_i^k $ we use  step by step triangles which are nearest $ x_i $ and we take the mediatrices of those triangles.

$$ \Omega_k=\{ x \in B(x_i, \delta_i\epsilon), |x-t_i^k| \leq |x-t_i^j|, j \not = k \}, $$

we write:

$$ B_i = \sum_{k=1}^m \int_{\Omega_k}  <(x-x_i) |\nabla u_i >(V_i(x)-V_i(x_i)) e^{u_i} dy. $$

We set,

$$ B_i^k= \int_{\Omega_k}  <(x-x_i) |\nabla u_i >(V_i(x)-V_i(x_i)) e^{u_i} dy, $$

We divide this integral in 4 integrals:

$$ B_i^k= \int_{\Omega_k}  <(x-t_i^k) |\nabla u_i >(V_i(x)-V_i(x_i)) e^{u_i} dy +\int_{\Omega_k}  <(t_i^k-x_i) |\nabla u_i >(V_i(x)-V_i(x_i)) e^{u_i} dy = $$

$$ = \int_{\Omega_k}  <(x-t_i^k) |\nabla u_i >(V_i(x)-V_i(t_i^k)) e^{u_i} dy + \int_{\Omega_k}  <(x-t_i^k) |\nabla u_i >(V_i(t_i^k)-V_i(x_i)) e^{u_i} dy + $$

$$ + \int_{\Omega_k}  <(t_i^k-x_i) |\nabla u_i >(V_i(x)-V_i(t_i^k)) e^{u_i} dy + \int_{\Omega_k}  <(t_i^k-x_i) |\nabla u_i >(V_i(t_i^k)-V_i(x_i)) e^{u_i} dy, $$

We set:

$$ A_1=\int_{\Omega_k}  <(x-t_i^k) |\nabla u_i >(V_i(x)-V_i(t_i^k)) e^{u_i} dy, $$

$$ A_2=\int_{\Omega_k}  <(x-t_i^k) |\nabla u_i >(V_i(t_i^k)-V_i(x_i)) e^{u_i} dy, $$

$$ A_3= \int_{\Omega_k}  <(t_i^k-x_i) |\nabla u_i >(V_i(x)-V_i(t_i^k)) e^{u_i} dy,  $$

$$ A_4= \int_{\Omega_k}  <(t_i^k-x_i) |\nabla u_i >(V_i(t_i^k)-V_i(x_i)) e^{u_i} dy. $$

For $ A_1 $ and $ A_2 $ we use the fact that in $ \Omega_k $ we have:

$$ u_i(x)+ 2 \log |x-t_i^k| \leq C,  $$

to conclude that for $ 0 < s \leq 1 $:

$$ A_1=A_2=o(1), $$

we have integrals of the form:

$$ A_1'=\int_{\Omega_k}  |\nabla u_i| e^{(1/2-s/2)u_i} dy = o(1), $$

and,

$$ A_2'=\int_{\Omega_k}  |\nabla u_i| e^{(1/2-s/4)u_i} dy= o(1). $$

For $ A_3 $ we use the previous fact and the $ \sup + \inf $ inequality to conclude that for $ 1/2 < s \leq 1 $:

$$ A_3= o(1) $$

because we have an integral of the form:

$$ A_3' = \int_{\Omega_k}  |\nabla u_i| e^{(3/4-s/2)u_i} dy = o(1). $$

For $ A_4 $ we use integration by part to have:

$$ A_4= \int_{\partial \Omega_k}  <(t_i^k-x_i) | \nu >(V_i(t_i^k)-V_i(x_i)) e^{u_i} dy. $$

But, the boundary of $ \Omega_k $ is the union of parts of mediatrices of segments linked to $ t_i^k $. Let's consider  a point $ t_i^j $ linked to  $ t_i^k $ and denote $ D_{i, j, k} $ the mediatrice of  the segment $ (t_i^j, t_i^k) $, which is in the boundary of $ \Omega_k $. Note that this mediatrice is in the boundary of $ \Omega_j $ and the same decompostion for $ \Omega_j $ gives us  the following term:

$$ A_4'= - \int_{D_{i, j, k}}  <(t_i^j-x_i) | \nu >(V_i(t_i^j)-V_i(x_i)) e^{u_i} dy. $$

Thus, we have  to estimate the sum of the 2 following terms:

$$ A_5=  \int_{D_{i, j, k}}  <(t_i^k-x_i) | \nu >(V_i(t_i^k)-V_i(x_i)) e^{u_i} dy. $$

and,

$$ A_6= A_4'= - \int_{D_{i, j, k}}  <(t_i^j-x_i) | \nu >(V_i(t_i^j)-V_i(x_i)) e^{u_i} dy. $$

We can write them as follows:

$$ A_5=  \int_{D_{i, j, k}}  <(x-x_i) | \nu >(V_i(t_i^k)-V_i(x_i)) e^{u_i} dy + \int_{D_{i, j, k}}  <(t_i^k-x) | \nu >(V_i(t_i^k)-V_i(x_i)) e^{u_i} dy. $$

and,

$$ A_6=  -\int_{D_{i, j, k}}  <(x-x_i) | \nu >(V_i(t_i^j)-V_i(x_i)) e^{u_i} dy -\int_{D_{i, j, k}}  <(t_i^j-x) | \nu >(V_i(t_i^j)-V_i(x_i)) e^{u_i} dy. $$

We can write:

$$ \int_{D_{i, j, k}}  <(x-x_i) | \nu >(V_i(t_i^k)-V_i(x_i)) e^{u_i} dy -  \int_{D_{i, j, k}}  <(x-x_i) | \nu >(V_i(t_i^j)-V_i(x_i)) e^{u_i} dy = $$

$$ =  \int_{D_{i, j, k}}  <(x-x_i) | \nu >(V_i(t_i^k)-V_i(x_i^j)) e^{u_i} dy = o(1), $$

for $ 1/2 < s \leq 1 $. Because, we do integration on the mediatrice of $ (t_i^j, t_i^k) $, $ |x-t_i^j|=|x-t_i^k| $, and:

$$ |V_i(t_i^k)-V_i(x_i^j)| \leq 2A |x-t_i^k|^s $$

$$ u_i(x)+ 2 \log |x-t_i^k| \leq C,  $$

and,

$$  |x- x_i| \leq \delta_i \epsilon,  $$

To estimate the integral of the  following term:

$$ e^{(3/4-s/2) u_i} \leq C r^{(-3/2+s)}, $$

which is intgrable and tends to 0, for $ 1/2 < s \leq 1 $, because we are on the ball $ B(x_i,  \delta_i  \epsilon) $.

In other part, for the term:

$$ \int_{D_{i, j, k}}  <(t_i^k-x) | \nu >(V_i(t_i^k)-V_i(x_i)) e^{u_i} dy- \int_{D_{i, j, k}}  <(t_i^j-x) | \nu >(V_i(t_i^j)-V_i(x_i)) e^{u_i} dy. $$

We use the fact that, on $ D_{i, j, k} $:

$$  |x-t_i^j|=|x-t_i^k|, $$

$$ u_i(x)+ 2 \log |x-t_i^k| \leq C,  $$

$$ |V_i(t_i^k)-V_i(x_i)| \leq 2A |x_i -t_i^k|^s \leq \delta_i^s, $$

and,
 
$$ |V_i(t_i^j)-V_i(x_i)| \leq 2A |x_i- t_i^j|^s \leq \delta_i^s, $$

To estimate the integral of the  following term:

$$ e^{(1/2-s/4) u_i} \leq C r^{(-1+s/2)}, $$

which is intgrable and tends to 0, because we are on the ball $ B(x_i,  \delta_i  \epsilon) $.

Thus,

$$ B_i=o(1), $$

\underbar {Proof of the theorem 2:}

\bigskip

Next, we use the formulation of the case of three blow-up points, see [2]. Because the blow-ups points are simple, we can consider the following function:

$$ v_i(\theta)= u_i( x_i+r_i\theta)-u_i(x_i),  $$
 
 where $ r_i $ is such that:
 
 $$ r_i = e^{-u_i(x_i)/2}, $$

$$ \int_{B^{\epsilon}} V_i(x_i+ \delta_i y) e^{v_i} \to 8 \pi. $$

$$ u_i(x_i+r_i \theta)=\int_{\Omega} G(x_i+r_i\theta,y) V_i(y) e^{u_i(y)} dx= $$

$$ =\int_{\Omega-B(x_i, 2\delta_i\epsilon')} G(x_i,y) V_i e^{u_i(y)} dy + \int_{B(x_i, 2\delta_i\epsilon')} G(x_i+r_i\theta,y) V_ie^{u_i(y)} dy = $$ 
 
We write, $ y=x_i+r_i\tilde \theta $, with $ |\tilde \theta| \leq 2\dfrac{\delta_i}{r_i}\epsilon' $,

 $$u_i(x_i+r_i \theta) =\int_{B(0, 2\frac{\delta_i}{r_i}\epsilon')} \dfrac{1}{2 \pi} \log \dfrac{|1-(\bar x_i+ r_i \bar \theta )(x_i+r_i \tilde \theta )|}{r_i|\theta-\tilde \theta|} V_ie^{u_i(y)} r_i^2dy  + $$
 
 $$ + \int_{\Omega-B(x_i, 2\delta_i\epsilon')} G(x_i+r_i\theta,y) V_i e^{u_i(y)} dy $$

 $$ u_i(x_i)=\int_{\Omega-B(x_i, 2\delta_i\epsilon')} G(x_i,y) V_i e^{u_i(y)} dy + \int_{B(x_i, 2\delta_i\epsilon' )} G(x_i,y) V_ie^{u_i(y)} dy $$
  
 Hence,
 
 $$  u_i(x_i)=\int_{B(0, 2\frac{\delta_i}{r_i}\epsilon')} \dfrac{1}{2 \pi} \log \dfrac{|1-\bar x_i(x_i+r_i \tilde \theta )|}{r_i|\tilde \theta|} V_ie^{u_i(y)} r_i^2dy  + $$
 
 $$ + \int_{\Omega-B(x_i,2\delta_i\epsilon')} G(x_i,y) V_i e^{u_i(y)} dy  $$
   
We look to the difference,

$$ v_i(\theta)= u_i( x_i+r_i\theta)-u_i(x_i)=\int_{B(0, 2\frac{\delta_i}{r_i}\epsilon')} \dfrac{1}{2 \pi} \log \dfrac{|\tilde \theta |}{|\theta-\tilde \theta|} V_ie^{u_i(y)} r_i^2dy  +  h_1+ h_2, $$

where,

 $$ h_1(\theta) = \int_{\Omega-B(x_i, 2\delta_i\epsilon')} G(x_i+r_i\theta,y) V_i e^{u_i(y)} dy -  \int_{\Omega-B(x_i, 2\delta_i\epsilon')} G(x_i,y) V_i e^{u_i(y)} dy, $$  
 
 and,
 
 $$ h_2(\theta)= \int_{B(0, 2\delta_i\epsilon')} \dfrac{1}{2 \pi} \log \dfrac{|1-(\bar x_i+ r_i \bar \theta )y|}{|1-\bar x_i y |} V_ie^{u_i(y)} dy. $$
 
 Remark that, $ h_1 $ and $ h_2 $ are two harmonic functions, uniformly bounded.
 
 \bigskip
  
According to the maximum principle,  the harmonic function $ G(x_i+r_i\theta, .) $ on $ \Omega-B(x_i,2\delta_i\epsilon') $ take its maximum on the boundary of $ B(x_i, 2\delta_i\epsilon') $, we can compute this maximum:

$$ G(x_i+r_i\theta, y_i)=\dfrac{1}{2 \pi} \log \dfrac{|1-(\bar x_i+ r_i\bar \theta)y_i|}{|x_i+r_i\theta-y_i|} \simeq \dfrac{1}{2\pi} \log \dfrac{(|1+|x_i|)\delta_i-\delta_i(3\epsilon'+o(1))|}{\delta_i \epsilon' } \leq C_{\epsilon'} < + \infty $$

with $ y_i=x_i+ 2\delta_i \theta_i \epsilon' $,  $ |\theta_i|= 1 $, and $ |r_i\theta| \leq \delta_i\epsilon' $.

\bigskip

We can remark, for $ |\theta| \leq \dfrac{\delta_i\epsilon'}{r_i}  $, that $ v_i $ is such that:

$$ v_i = h_1+h_2+\int_{B(0, 2\frac{\delta_i}{r_i}\epsilon')} \dfrac{1}{2 \pi} \log \dfrac{|\tilde \theta |}{|\theta-\tilde \theta|} V_ie^{u_i(y)} r_i^2dy, $$

$$ v_i = h_1+h_2+\int_{B(0,2\frac{\delta_i}{r_i} \epsilon')} \dfrac{1}{2 \pi} \log \dfrac{|\tilde \theta |}{|\theta-\tilde \theta|} V_i (x_i+ r_i \tilde \theta) e^{v_i(\tilde \theta)} d\tilde \theta, $$

with $ h_1 $ and $ h_2 $, the two uniformly bounded harmonic functions.

\bigskip

{\bf Remark:} In the case of 2 or 3 or 4 blow-up points, and if we consider the half ball, we have supplemntary terms,  around the 2 other blow-up terms. Note that the Green function of the half ball is quasi-similar to the one of the unit ball and our computations are the same if we consider the half ball.

\bigskip

By the asymptotic estimates of Cheng-Lin, we can see that, we have the following uniform estimates at infinity. We have, after considering the half ball and its Green function, the following estimates:

$ \forall \,\, \epsilon>0, \,\epsilon'>0 \,\, \exists \,\, k_{\epsilon,\epsilon'} \in {\mathbb R}_+ , \,\, i_{\epsilon, \epsilon'} \in {\mathbb N} $ and $ C_{\epsilon, \epsilon'} >0 $, such that, for $ i \geq i_{\epsilon,\epsilon'} $ and $ k_{\epsilon, \epsilon'} \leq  |\theta| \leq \dfrac{\delta_i\epsilon'}{r_i}  $,
 
 $$  (-4-\epsilon)\log | \theta | -C_{\epsilon, \epsilon'}  \leq v_i(\theta) \leq (-4+\epsilon)\log | \theta | + C_{\epsilon, \epsilon'}, $$
 
 and,
 
$$ \partial_j v_i \simeq \partial_j u_0(\theta) \pm \dfrac{\epsilon}{|\theta|} + C\left (\dfrac{r_i}{\delta_i} \right )^2|\theta| + m\times \left (\dfrac{r_i}{\delta_i} \right ) +  $$

$$ + \sum_{k=2}^m  C_1\left (\dfrac{r_i}{d(x_i, x_i^k)} \right ), $$

In the case, we have:

$$  \dfrac{d(x_i,x_i^k)}{\delta_i}\to + \infty \,\,\, {\rm for } \,\,\, k =2\ldots m, $$

We have after using the  previous term of the Pohozaev identity, for $ 1/2 < s \leq 1 $:
 
 $$ o(1)=J'_i= m'+ \sum_{k=1}^m C_k o(1), $$
 
 $$ 0=  \lim_{\epsilon' } \lim_{\epsilon} \lim_i J_i'= m', $$
 
 which contradict the fact that $ m'>0 $.

here, 

$$  J_i= B_i= \int_{B(x_i, \delta_i \epsilon')} <x_1^i |\nabla (u_i-u)> (V_i-V_i(x_i)) e^{u_i} dy. $$

We use the previous formulation around each blow-up point.

If, for $ x_i^j $, we have:

$$  \dfrac{d(x_i^j,x_i^k)}{\delta_i^j}\to + \infty \,\,\, {\rm for } \,\,\, k \not =j, k=1\ldots m, $$

We use the previous formulation around this blow-up point. We consider the following quantity:

$$  J_i^j= B_i^j= \int_{B(x_i^j, \delta_i^j \epsilon')} <x_1^{i,j} |\nabla (u_i-u)> (V_i-V_i(x_i^j)) e^{u_i} dy. $$

with, 

$$ x_1^{i,j} = (\delta_i^j, 0), $$

In this case, we set:

$$ v_i^j(\theta)= u_i( x_i^j+r_i^j\theta)-u_i(x_i^j),  $$
 
 where $ r_i^j $ is such that:
 
 $$ r_i^j = e^{-u_i(x_i^j)/2}, $$

$$ \int_{B(x_i^j, \delta_i^j \epsilon')} V_i(x_i^j+ \delta_i^j y) e^{v_i} \to 8 \pi. $$

We have, after considering the half ball and its Green function, the following estimates:

$ \forall \,\, \epsilon>0, \,\epsilon'>0 \,\, \exists \,\, k_{\epsilon,\epsilon'} \in {\mathbb R}_+ , \,\, i_{\epsilon, \epsilon'} \in {\mathbb N} $ and $ C_{\epsilon, \epsilon'} >0 $, such that, for $ i \geq i_{\epsilon,\epsilon'} $ and $ k_{\epsilon, \epsilon'} \leq  |\theta| \leq \dfrac{\delta_i^j\epsilon'}{r_i^j}  $,
 
 $$  (-4-\epsilon)\log | \theta | -C_{\epsilon, \epsilon'}  \leq v_i^j(\theta) \leq (-4+\epsilon)\log | \theta | + C_{\epsilon, \epsilon'}, $$
 
 and,
 
$$ \partial_k v_i^j \simeq \partial_k u_0^j(\theta) \pm \dfrac{\epsilon}{|\theta|} + C\left (\dfrac{r_i^j}{\delta_i^j} \right )^2|\theta| + m\times \left (\dfrac{r_i^j}{\delta_i^j} \right ) +  $$

$$ + \sum_{l \not = j}^m  C_1\left (\dfrac{r_i^j}{d(x_i^j, x_i^l)} \right ), $$

We have after using the  previous term of the Pohozaev identity, for $ 1/2 < s \leq 1 $:
 
 $$ o(1)=J_i^j= B_i^j= m'+ \sum_{ l \not =j}^m C_l o(1), $$
 
 $$ 0=  \lim_{\epsilon' } \lim_{\epsilon} \lim_i J_i^j= m', $$
 
 which contradict the fact that $ m'>0 $.

\bigskip

If, for $ x_i^j $, we have:

$$  \dfrac{d(x_i^j,x_i^k)}{\delta_i^j}\leq C_{j, k} \,\,\, {\rm for \,\, some } \,\,\, k=k_j \not =j,  1\leq k \leq m,$$

All the distances $  d(x_i^j,x_i^k) $ are comparable with some  $ \delta_i^j $. This means that we can use the Pohozaev identity directly. We can do this for example, for 4 blow-ups points.

\bigskip

We have many cases:

\smallskip

Case 1: the blow-up points are "equivalents", it seems that we have the same radius for the blow-up points.

\smallskip

Case 2: 3 points are "equivalents" and another blow-up point linked to the 3 blow-up points. We apply the Pohozaev identity directly with central point which link the 3 blow-up to the last.

\smallskip

Case 3: 2 pair of blow-up points separated.

\smallskip

Case 3.1: the 2 pair are linked: we apply the Pohozaev identity.

\smallskip

Case 3.2: the two pair are separated. It is the case of two separated blow-up points, see [1]

\bigskip 

\begin{center}

{\bf  ACKNOWLEDGEMENT. }

\end{center}

\smallskip

The author is grateful to Professor G. Wolansky who has communicated him the private communcation.



 
 







\end{document}